\newif\ifarxiv

\arxivtrue

\ifarxiv
	\RequirePackage{fix-cm}
	\documentclass[smallextended]{svjour3}
	\usepackage{fullpage}
	\smartqed
\else
	\documentclass[final,5p,times,twocolumn]{elsarticle}
\fi

\usepackage{amssymb}

\usepackage{xcolor,hyperref}
\hypersetup{
        linktoc=page,
        colorlinks=true,
        linkcolor=blue,
        citecolor=red,
        filecolor=magenta,
        urlcolor=cyan
}

\usepackage{amsmath}

\ifarxiv
\else
        \journal{Operations Research Letters}

	\usepackage{amsthm}
        \newtheorem{theorem}{Theorem}
        \newtheorem{lemma}{Lemma}
        \newtheorem{corollary}{Corollary}
        
        \newtheorem{example}{Example}
        \newtheorem{remark}{Remark}
\fi
\newtheorem{assumption}{Assumption}

\numberwithin{equation}{section}
\numberwithin{theorem}{section}
\numberwithin{lemma}{section}
\numberwithin{corollary}{section}
\numberwithin{definition}{section}
\numberwithin{example}{section}
\numberwithin{remark}{section}

\newcommand{\dom} {\operatorname{dom}}
\newcommand{\prox} {\operatorname{prox}}

\newcommand{\pL} {\operatorname{p}_L^{ }}
\newcommand{\pLk} {\operatorname{p}_{L_k}^{ }}
\newcommand{\pLi} {\operatorname{p}_{L_i}^{ }}
\newcommand{\pLkp} {\operatorname{p}_{L_{k}'}^{ }}

\newcommand{\E} {\mathbb{R}^n} 
\newcommand{\V} {\mathbb{R}^m} 

\newcommand{\x} {\mathbf{x}}
\newcommand{\y} {\mathbf{y}}
\newcommand{\z} {\mathbf{z}}
\newcommand{\w} {\mathbf{w}}
\renewcommand{\u} {\mathbf{u}}
\renewcommand{\v} {\mathbf{v}}
\renewcommand{\b} {\mathbf{b}}
\renewcommand{\d} {\mathbf{d}}
\newcommand{\s} {\mathbf{s}}
\newcommand{\bz} {\bar{\z}}

\newcommand{\tH} {\tilde{H}}
\newcommand{\tq} {\tilde{q}}

\newcommand{\Aop} {\mathbf{A}} 
\newcommand{\X} {\mathcal{X}}

\DeclareMathOperator*{\argmin}{argmin}
\DeclareMathOperator*{\argmax}{argmax}

\newcommand{\TheTitle}{Fast dual proximal gradient algorithms
with rate $O(1/k^{1.5})$ for convex minimization}

\begin{document}

\ifarxiv
	\title{\TheTitle\thanks{This research was supported in part by NIH grant U01 EB018753.}}
	\author{Donghwan Kim \and Jeffrey A. Fessler}
        \institute{Donghwan Kim \and Jeffrey A. Fessler \at
              Dept. of Electrical Engineering and Computer Science,
                University of Michigan, Ann Arbor, MI 48109 USA \\
              \email{kimdongh@umich.edu, fessler@umich.edu}          
        }
        \date{Date of current version: \today}

	\maketitle
\else
	\begin{frontmatter}

	\title{\TheTitle}

	\author{Donghwan Kim\corref{cor1}}
	\cortext[cor1]{Corresponding author}
	\ead{kimdongh@umich.edu}
	
	\author{Jeffrey A. Fessler\corref{}}
	\ead{fessler@umich.edu}

	\address{Department of Electrical Engineering and Computer Science,
	University of Michigan, Ann Arbor, MI, USA}
\fi

\begin{abstract}
%
We consider minimizing the composite function 
that consists of a strongly convex function and a convex function. 
The fast dual proximal gradient (FDPG) method decreases 
the dual function with a rate $O(1/k^2)$, 
leading to a rate $O(1/k)$ for decreasing the primal function. 
We propose a generalized FDPG method that guarantees an $O(1/k^{1.5})$ rate 
for the dual proximal gradient norm decrease. 
By relating this to the primal function decrease, 
the proposed approach decreases the primal function with the improved $O(1/k^{1.5})$ rate.
\end{abstract}

\ifarxiv
\else
	\end{frontmatter}
\fi


\section{Introduction}

This paper focuses on improving the rate of convergence
of dual-based proximal gradient methods
for minimizing the sum of two convex functions,
where one is assumed to be strongly convex.
The convergence analysis in this paper
focuses on the rate of decrease of the dual proximal gradient norm,
whereas the existing analysis in~\cite{beck:14:afd} 
focuses on the rate of decrease of the dual function.

This work is based on
the alternating minimization algorithm by Tseng~\cite{tseng:91:aoa}
that exploits the strong convexity.
The method~\cite{tseng:91:aoa}
is essentially equivalent to applying the proximal gradient method
to the dual function,
which is naturally named a dual proximal gradient (DPG) method 
in~\cite{beck:14:afd}.
In~\cite{beck:14:afd,goldstein:14:fad},
this alternating minimization algorithm (or DPG) 
is accelerated using the fast proximal gradient method (FPGM) in~\cite{beck:09:afi},
widely popularized under the name FISTA.
That fast DPG (FDPG) method
decreases the dual function at rate $O(1/k^2)$
due to the acceleration of FISTA~\cite{beck:09:afi},
where $k$ denotes the number of iterations;
the FDPG method is effective for various applications
such as total-variation-based image denoising problems
\cite{goldstein:14:fad,beck:09:fgb}
and model predictive control problems
\cite{pu:14:fam}.

In the interest of the \emph{primal} convergence analysis of
DPG and FDPG methods,
Beck and Teboulle~\cite{beck:14:afd} derived
nonasymptotic convergence bounds
for the decrease of the distance
between the primal sequence and a primal solution,
and for the primal function decrease of DPG and FDPG.
In particular,
the rate $O(1/k^2)$ for the dual function decrease of FDPG
provided the rate $O(1/k)$ for both the primal distance and function decrease,
which is superior to those rates of subgradient and DPG methods
in~\cite{beck:14:afd}.

In addition to analyzing the primal convergence analysis using 
the dual function decrease as in~\cite{beck:14:afd},
Nesterov~\cite{nesterov:12:htm} 
pointed out that the dual gradient decrease
is closely related to the primal function decrease
for minimizing a strongly convex function
with a linear equality constraint.
He then suggested using an algorithm
that decreases the dual gradient with a fast rate $O(1/k^{1.5})$,
thus providing the same rate for the primal function decrease. 
That analysis was extended to a linear inequality constrained strongly convex problem
in~\cite{necoara:16:ica}.
This paper further extends such analyses to strongly convex composite problems,
by showing that the dual \emph{proximal} gradient decrease
is directly related to the primal function decrease.

We recently proposed an accelerated proximal gradient method 
named the generalized FPGM (GFPGM) in~\cite{kim:16:ala}
that has rate $O(1/k^{1.5})$
for decreasing the proximal gradient norm
and that is computationally as efficient as FISTA.
This paper proposes to incorporate that method with duality,
leading to a generalized FDPG (GFDPG) method.
We show that the proposed approach 
has the rate $O(1/k^{1.5})$ for decreasing the primal function,
by extending the analysis in~\cite{nesterov:12:htm,necoara:16:ica}.
As a byproduct of our analysis,
we prove an $O(1/k)$ bound 
on the rate of decrease of the primal function,
which is interestingly the same as that of the FDPG in~\cite{beck:14:afd}.

Sometimes the function information such as the strong convexity parameter is unavailable
or difficult to approximate,
and the FDPG method (and FISTA) have a backtracking scheme~\cite{beck:09:afi}
that circumvents that problem.
By introducing such a backtracking scheme to GFPGM~\cite{kim:16:ala},
we illustrate that the proposed GFDPG
also satisfies an $O(1/k^{1.5})$ bound 
on the primal function decrease for such cases.

Section~\ref{sec:opt} presents 
the optimization problem of interest and its dual.
Section~\ref{sec:fdpg} reviews 
the convergence analysis of FDPG in~\cite{beck:14:afd}.
Section~\ref{sec:rate} analyzes the convergence rate
of the primal function decrease
using the dual proximal gradient norm convergence.
Section~\ref{sec:gfdpg} proposes using
the accelerated proximal gradient method named GFPGM in~\cite{kim:16:ala}
instead of FISTA to effectively tackle the dual problem,
leading to an improved $O(1/k^{1.5})$ rate
for the primal function decrease.
Section~\ref{sec:conc} concludes.


\section{Optimization problem and its dual}
\label{sec:opt}

\subsection{The problem}

This paper considers the following composite convex problem:
\begin{align}
\x_* := \argmin_{\x} \left\{H(\x) := f(\x) + g(\Aop\x)\right\}
\tag{P}
\label{eq:P}
,\end{align}
where both $f\,:\,\E\to(-\infty,+\infty]$ and $g\,:\,\V\to(-\infty,+\infty]$ 
are proper, closed, and convex extended real-valued functions,
while the function $f$ is further assumed to be $\sigma$-strongly convex for $\sigma>0$,
and $\Aop$ is a $m\times n$ matrix.
Due to the strong convexity, problem~\eqref{eq:P}
has a unique optimal solution $\x_*$.

Problem~\eqref{eq:P} is general enough
to model various applications;
representative examples
such as image denoising, projection onto the intersection of convex sets,
and resource allocation problems
are provided in~\cite{beck:14:afd}
(see also~\cite{goldstein:14:fad,beck:09:fgb,pu:14:fam}).
Tackling such problems directly (in a primal domain)
using algorithms such as subgradient methods
suffer from relatively slow convergence rates~\cite{beck:14:afd}.
The next subsection and Section~\ref{sec:fdpg}
review the fast proximal gradient scheme combined with duality
in~\cite{beck:14:afd,goldstein:14:fad}
that exploits the properties of problem~\eqref{eq:P}
and that converges faster than the subgradient methods~\cite{beck:14:afd}.

\subsection{The dual problem}

Problem~\eqref{eq:P}
has the following equivalent constrained form:
\begin{align}
\x_* = \argmin_{\x}\min_{\z} \left\{\tH(\x,\z) := f(\x) + g(\z) \,:\, \Aop\x - \z = 0\right\}
\tag{P$'$}
\label{eq:P_}
,\end{align}
where $H(\x) = \tH(\x,\Aop\x)$.
Problem~\eqref{eq:P_} has the following dual problem:
\begin{align}
\y_* \in \argmax_{\y} q(\y)
\tag{D}
\label{eq:D}
,\end{align}
where the dual function is defined as
\cite{beck:14:afd}:
\begin{align}
q(\y) &:= \min_{\x,\z} \left\{\tH(\x,\z) - \langle\y,\Aop\x - \z\rangle\right\}
= -f^*(\Aop^\top\y) - g^*(-\y)
\end{align}
with dual variable vector $\y\in\V$.
Let $\y_*$ denote an optimal (dual) solution of problem~\eqref{eq:D}.
The convex conjugates of $f$ and $g$ are defined as
\begin{align*}
f^*(\u) = \max_{\x} \{\langle\u,\x\rangle - f(\x)\},
\;\text{and}\;
g^*(\y) = \max_{\z} \{\langle\y,\z\rangle - g(\z)\}
.\end{align*}

To make problem~\eqref{eq:D} into an equivalent convex problem for convenience,
\cite{beck:14:afd} defines 
\begin{align}
F(\y) := f^*(\Aop^\top\y),
\;\text{and}\;
G(\y) := g^*(-\y)
,\end{align}
where $F$ has a Lipschitz continuous gradient
(due to the strong convexity of $f$)
with a constant $L_F := \frac{||\Aop||^2}{\sigma}$
\cite[Lemma 3.1]{beck:14:afd}, {\it i.e.},
for any $\x,\u\in\E$
\begin{align}
||\nabla F(\x) - \nabla F(\u)|| \le L_F||\x - \u||
.\end{align}
Then dual problem~\eqref{eq:D} is equivalent to the following:
\begin{align}
\y_* \in \argmin_{\y} \left\{\tq(\y) := F(\y) + G(\y)\right\}
\tag{D$'$}
\label{eq:D_}
\end{align}
that consists of a smooth function $F$
and a closed proper function $G$.
One can solve using proximal gradient methods.
Note that $\tq(\y) = - q(\y)$ by definition.

Even when solving the dual problem~\eqref{eq:D} (or~\eqref{eq:D_}),
we are eventually interested in analyzing 
the convergence rate of the primal sequence
as in~\cite{beck:14:afd} and this paper.
For a given dual variables vector $\y$,
the corresponding primal variables vectors
are defined as
$(\x(\y),\z(\y)) \in 
\argmin_{\x,\z} \left\{\tH(\x,\z) - \langle\y,\Aop\x - \z\rangle\right\}$,
{\it i.e.},
\begin{align}
\x(\y) &:= \argmax_\x \left\{\langle\Aop^\top\y,\x\rangle - f(\x)\right\}, 
	\label{eq:primx} \\
\z(\y) &\in \argmin_\z \left\{\langle\y,\z\rangle + g(\z)\right\} 
	\label{eq:primz}
.\end{align}
Then by definition, $\x_* = \x(\y_*)$ and these vectors satisfy
\begin{align}
\tH(\x(\y),\z(\y)) - q(\y) = \langle\y,\Aop\x(\y) - \z(\y)\rangle
\label{eq:pd}
.\end{align}

Next, Section~\ref{sec:fdpg} reviews
bounds on the convergence rate of the primal function decrease 
for the primal variable vector $\x(\y)$
of dual-based proximal gradient methods
using bounds on the dual function decrease~\cite{beck:14:afd}.
In contrast, Sections~\ref{sec:rate} and~\ref{sec:gfdpg}
analyze the primal sequence
using~\eqref{eq:pd}
and bounds on the dual \emph{proximal} gradient decrease.

\section{Fast dual-based proximal gradient methods}
\label{sec:fdpg}

\subsection{Dual-based proximal gradient methods}

The proximal gradient method~\cite{beck:09:afi} 
for solving~\eqref{eq:D_}
has the following update at $k$th iteration for $k\ge1$
with given $L_0$ and $\y_0$:\footnote{
The Moreau proximal map~\cite{moreau:65:ped}
of a proper closed and convex function $h\,:\,\V\to(-\infty,\infty]$
in~\eqref{eq:proxgrad} is defined as
$
\prox_h(\w) = \argmin_{\y\in\V} \left\{h(\y) + \frac{1}{2}||\y - \w||^2\right\}.
$
}
\begin{align}
\y_k &= \pLk(\y_{k-1}) \nonumber \\ 
	&:= \argmin_\y \left\{
	\begin{array}{l}
		Q_{L_k}(\y,\y_{k-1}) :=
		F(\y_{k-1}) + \langle\y-\y_{k-1},\nabla F(\y_{k-1})\rangle \\ 
		\qquad\qquad\qquad\qquad\quad\;+ \frac{L_k}{2}||\y-\y_{k-1}||^2 + G(\y)
	\end{array}\right\} \nonumber \\
	&=\prox_{\frac{1}{L_k}G} \left(\y_{k-1} - \frac{1}{L_k}\nabla F(\y_{k-1})\right)
\label{eq:proxgrad}
,\end{align}
where $L_k$ is chosen to satisfy $L_{k-1} \le L_k$
and $\tq(\pLk(\y_{k-1})) \le Q_{L_k}(\pLk(\y_{k-1}),\y_{k-1})$,
which guarantees descent because 
$Q_{L_k}(\pLk(\y_{k-1}),\y_{k-1}) \le Q_{L_k}(\y_{k-1},\y_{k-1}) = \tq(\y_{k-1})$.
Using the fixed constant $L_k = L_F$ for all $k$
can satisfy the condition on $L_k$.
However when $L_F$ is unknown
or cannot be easily approximated,
a backtracking scheme in~\cite{beck:09:afi}
can be adopted.
This proximal gradient method
decreases the (dual) function with rate $O(1/k)$~\cite{beck:09:afi}.

The proximal gradient update $\pLk(\y_{k-1})$ in~\eqref{eq:proxgrad}
has an equivalent efficient update
in terms of the original functions $f$ and $g$ as follows~\cite[Lemma 3.2]{beck:14:afd}: 
\begin{align}
\u_k &= \x(\y_{k-1}), \label{eq:uk} \\ 
\v_k &= \prox_{L_kg} \left(\Aop\u_k - L_k\y_{k-1}\right), \label{eq:vk} \\
\y_k &= \y_{k-1} - \frac{1}{L_k}\left(\Aop\u_k - \v_k\right) \label{eq:yk}
,\end{align}
which exactly matches the update of the alternating minimization algorithm 
in~\cite{tseng:91:aoa}.
The advantage of this alternating minimization algorithm
over the augmented Lagrangian-based methods~\cite{glowinski:89}
for solving~\eqref{eq:P} (or~\eqref{eq:P_})
is that the method can exploit separability of $f$
in the update step~\eqref{eq:uk}.

The next section reviews FDPG~\cite{beck:14:afd,goldstein:14:fad}, 
the accelerated version of DPG using FISTA~\cite{beck:09:afi}.

\subsection{FDPG method and its convergence analysis}

In~\cite{beck:14:afd,goldstein:14:fad},
DPG is accelerated using FISTA~\cite{beck:09:afi}
with negligible extra computation per iteration as shown below,
which is named FDPG.

\begin{table}[h]
\ifarxiv
	\normalsize
\fi
\centering
\begin{tabular}{|l|}
\hline
The FDPG Method with backtracking \\
\hline
Input: Take $L_0$, $\y_0 = \w_0$, $t_0 = 1$. \\
Step $k$. ($k \ge 1$) \\
\quad Choose $L_k$ s.t. $L_{k-1} \le L_k$, and \\
\quad \qquad\qquad $\tq(\pLk(\w_{k-1})) \le Q_{L_k}(\pLk(\w_{k-1}),\w_{k-1})$. \\
\quad $\y_k = \pLk(\w_{k-1})$ \\
\quad $t_k = \frac{1 + \sqrt{1 + 4t_{k-1}^2}}{2}$ \\
\quad $\w_k = \y_k + \frac{t_{k-1} - 1}{t_k}(\y_k - \y_{k-1})$ \\
\hline
\end{tabular}
\end{table}

This FDPG has the following bound on the dual function decrease
with rate $O(1/k^2)$~\cite[Theorem 4.4]{beck:09:afi}, {\it i.e.},
\begin{align}
q(\y_*) - q(\y_k) = \tq(\y_k) - \tq(\y_*)
\le \frac{2L_k||\y_0 - \y_*||^2}{t_{k-1}^2}
\le \frac{2L_k||\y_0 - \y_*||^2}{(k+1)^2}
\label{eq:dualcost_fdpg}
.\end{align}
This rate is superior to the rate $O(1/k)$ 
for the dual function decrease of DPG~\cite[Theorem 3.1]{beck:09:afi}.


In~\cite{beck:14:afd},
it is shown that the rate $O(1/k^2)$ of the dual function decrease in~\eqref{eq:dualcost_fdpg}
provides the $O(1/k)$ bound on the convergence of the primal distance and function decrease.
%
In particular,
with the following assumption:\footnote{
\cite{beck:14:afd}
defines the closed and convex feasibility set
$\X = \{\x\in\E\,:\,\x\in\dom(f),\,\Aop\x\in\dom(g)\}$
and assumes $\gamma_H^{ } := \max_{\x\in\X}\max_{\d\in\partial H(\x)} ||\d|| < \infty$,
whereas this paper uses $\X = \E$.
}
\begin{assumption}
\label{assump}
The function $H$ is subdifferentiable for all $\x\in\E$,
and its subgradients are bounded as
\begin{align*}
\gamma_H^{ } := \max_{\x\in\E}\max_{\d\in\partial H(\x)} ||\d|| < \infty,
\end{align*}
\end{assumption}
\noindent
the corresponding primal sequence $\{\x(\y_k)\}$ of FDPG
defined by~\eqref{eq:primx} decreases
the primal function with rate $O(1/k)$
\cite[Theorem 4.3]{beck:14:afd}, {\it i.e.},
\begin{align}
H(\x(\y_k)) - H(\x_*) \le 2\gamma_H^{ }\sqrt{\frac{L_k}{\sigma}}\frac{||\y_0 - \y_*||}{k+1}
\label{eq:assump_cost_fdpg}
.\end{align}
%
%
In addition, the proof of~\cite[Theorem 4.3]{beck:14:afd} for~\eqref{eq:assump_cost_fdpg}
implies the following $O(1/\!\sqrt{k})$ bound
for the primal function decrease of DPG.
\begin{theorem}
\label{thm:dpg}
Let $\{\y_k\}$ be the sequence generated by DPG.
Then for any $k\ge1$ and with Assumption~\ref{assump},
the corresponding primal sequence $\{\x(\y_k)\}$ defined by~\eqref{eq:primx} satisfies
\begin{align}
H(\x(\y_k)) - H(\x_*) \le \gamma_H^{ }\sqrt{\frac{L_k}{\sigma}}\frac{||\y_0 - \y_*||}{\sqrt{k}}
\label{eq:assump_cost_dpg}
.\end{align}
\begin{proof}
This can be easily proven using~\cite[Theorem 3.1]{beck:09:afi}
that shows the $O(1/k)$ rate for the dual function decrease of DPG,
and using the proof of \cite[Theorems 4.1 and 4.3]{beck:14:afd}.
\ifarxiv
        \qed
\fi
\end{proof}
\end{theorem}

Both the bounds~\eqref{eq:assump_cost_fdpg} and~\eqref{eq:assump_cost_dpg}
resulting from the bound on the dual function decrease of FDPG and DPG respectively
seem to suggest that the primal function decrease of FDPG 
is faster than that of DPG.
However, the next section 
improves on~\eqref{eq:assump_cost_dpg}
by deriving an $O(1/k)$ bound
on the primal function decrease for DPG,
which is the same rate as that of FDPG in~\eqref{eq:assump_cost_fdpg}.
This new analysis in Section~\ref{sec:rate} 
uses a bound on the dual proximal gradient norm decrease
with an assumption that is weaker than Assumption~\ref{assump}
to analyze the primal function decrease.


\section{Rate of convergence of the primal function}
\label{sec:rate}

\subsection{Preliminaries}

This section presents two Lemmas that are the ingredients
for relating the dual proximal gradient norm $||\pL(\y) - \y||$
to the primal-dual gap $H(\x(\pL(\y))) - q(\pL(\y))$.
This in turn
determines the rate of the decrease
of the primal function $H(\x(\pL(\y))) - H(\x_*)$,
because $H(\x_*) = \tH(\x_*,\Aop\x_*) \ge q(\y_*) \ge q(\pL(\y))$.

\begin{lemma}
\label{lem:strconvex}
For any $\y,\w\in\V$, the following inequality holds:
\begin{align}
||\x(\y) - \x(\w)|| \le \frac{||\Aop||}{\sigma}||\y-\w||
.\end{align}
\begin{proof}
Since $f$ is $\sigma$-strongly convex, for any $\x,\u\in\E$ we have
\begin{align}
\sigma||\x - \u||
        &\le ||f'(\x) - f'(\u)|| \nonumber 
,\end{align}
where $f'(\x)\in\partial f(\x)$.
Then, using $\Aop^\top\y\in\partial f(\x(\y))$
that follows from the optimality condition of~\eqref{eq:primx},
we have
\begin{align}
\sigma||\x(\y) - \x(\w)||
        &\le ||\Aop^\top(\y - \w)|| \le  ||\Aop||\cdot||\y - \w|| \nonumber
.\end{align}
\ifarxiv
        \qed
\fi
\end{proof}
\end{lemma}

\begin{lemma}
\label{lem:proxcond}
For any $\y\in\V$ and $L>0$, the following equality holds:
\begin{align}
\Aop\x(\y) - \z(\pL(\y)) + L(\pL(\y) - \y) = 0
.\end{align}
\begin{proof}
%
%
We show that the following vector $\bz$: 
\begin{align}
\bz := L(\pL(\y) - \y) + \Aop\x(\y)
\label{eq:bz}
\end{align}
corresponds to $\z(\pL(\y))$.

Using~\eqref{eq:uk},~\eqref{eq:vk} and~\eqref{eq:yk}, we have
\begin{align}
\bz
= &\prox_{Lg} \left(\Aop\x(\y) - L\y\right) \nonumber \\
= &\argmin_\z \left\{Lg(\z) + \frac{1}{2}||\z - (\Aop\x(\y) - L\y)||^2\right\}
	\label{eq:bz_}
.\end{align}
The optimality condition of~\eqref{eq:bz_} implies
that there exists $g'(\bz) \in \partial g(\bz)$ such that
$
Lg'(\bz) + \bz - \Aop\x(\y) + L\y = 0
$
that is equivalent to
\begin{align}
g'(\bz) + \pL(\y) = 0
\label{eq:bz_cond}
\end{align}
using~\eqref{eq:bz}.
This condition~\eqref{eq:bz_cond} holds for $\bz = \z(\pL(\y))$
based on the optimality condition of
\begin{align*}
\z(\pL(\y)) \in \argmin_\z \{\langle\pL(\y),\z\rangle + g(\z)\}
\end{align*}
in~\eqref{eq:primz}, which concludes the proof.
\ifarxiv
        \qed
\fi
\end{proof}
\end{lemma}

\subsection{Relating the dual proximal gradient norm
to the primal-dual gap}

Based on Lemmas~\ref{lem:strconvex} and~\ref{lem:proxcond},
the following Lemma analyzes the convergence bound for
the primal-dual gap decrease 
$\tH(\x(\pL(\y)),\z(\pL(\y))) - q(\pL(\y))$ of~\eqref{eq:P_}.

\begin{lemma}
\label{lem:p_d_cost}
For any $\y\in\V$, $L>0$ and 
the corresponding primal vectors defined by~\eqref{eq:primx} and~\eqref{eq:primz}, 
the following inequality holds:
\begin{align}
&\tH(\x(\pL(\y)),\z(\pL(\y))) - q(\pL(\y)) \nonumber \\
&\qquad\qquad\qquad\qquad \le (L + L_F)\,||\pL(\y)||\cdot||\pL(\y) - \y||
.\end{align}
\begin{proof}
We have
\begin{align}
&\tH(\x(\pL(\y)),\z(\pL(\y))) - q(\pL(\y)) \nonumber \\
=\;& \langle\pL(\y),\,\Aop\x(\pL(\y)) - \z(\pL(\y))\rangle \nonumber \\
=\;& \langle\pL(\y),\,\Aop\x(\y) - \z(\pL(\y)) + \Aop(\x(\pL(\y)) - \x(\y))\rangle \nonumber \\
\le\;& ||\pL(\y)||\,(||\Aop\x(\y) - \z(\pL(\y))||
        + ||\Aop||\cdot||\x(\pL(\y)) - \x(\y)||) \nonumber \\
\le\;& ||\pL(\y)||\,\left(||L(\pL(\y) - \y)||
        + \frac{||\Aop||^2}{\sigma}||\pL(\y) - \y||\right) \nonumber \\
\le\;& ||\pL(\y)||
        \left(L + L_F\right)||\pL(\y) - \y||
\label{eq:pd_bound}
,\end{align}
where the first equality uses~\eqref{eq:pd},
the first inequality uses the Cauchy-Schwartz and triangle inequalities,
and the second inequality uses Lemmas~\ref{lem:strconvex} and~\ref{lem:proxcond}.
\ifarxiv
        \qed
\fi
\end{proof}
\end{lemma}

Lemma~\ref{lem:p_d_cost} shows that
the primal-dual gap decrease of~\eqref{eq:P_}
depends on the decrease of the dual proximal update $||\pL(\y) - \y||$.
However, we are more interested in the primal-dual gap decrease of~\eqref{eq:P}
than of~\eqref{eq:P_}.
Towards that end,
we introduce the following assumption 
that is weaker than Assumption~\ref{assump}.

\begin{assumption}
\label{assump_g}
The function $g$ is subdifferentiable for all $\z\in\V$,
and its subgradients are bounded as
\begin{align*}
\gamma_g := \max_{\z\in\V}\max_{\d\in\partial g(\z)} ||\d|| < \infty.
\end{align*}
\end{assumption}

We next analyze the convergence bound of
the primal-dual gap
$H(\x(\pL(\y))) - q(\pL(\y))$ of~\eqref{eq:P}
using Assumption~\ref{assump_g},
which is one of the main contribution of this paper.

\begin{lemma}
\label{lem:pd_cost}
With Assumption~\ref{assump_g},
for any $\y\in\V$, $L>0$ and 
the corresponding primal vector defined by~\eqref{eq:primx},
the following primal-dual gap inequality holds:
\begin{align}
&H(\x(\pL(\y))) - q(\pL(\y)) \nonumber \\
&\qquad\qquad\qquad \le (L + L_F)\,(||\pL(\y)|| + \gamma_g)\,||\pL(\y) - \y||
.\end{align}
\begin{proof}
We have
\begin{align}
&H(\x(\pL(\y))) - q(\pL(\y)) \nonumber \\
=\;& f(\x(\pL(\y))) + g(\Aop\x(\pL(\y))) - q(\pL(\y)) \nonumber \\
\le\;& f(\x(\pL(\y))) + g(\z(\pL(\y))) \nonumber \\
	&- \langle g'(\Aop\x(\pL(\y))),\,\z(\pL(\y)) - \Aop\x(\pL(\y))\rangle
        - q(\pL(\y)) \nonumber \\
=\;& \tH(\x(\pL(\y)),\z(\pL(\y))) - q(\pL(\y)) \nonumber \\
	&- \langle g'(\Aop\x(\pL(\y))),\,\z(\pL(\y)) - \Aop\x(\pL(\y))\rangle \nonumber \\
=\;& \langle\pL(\y) + g'(\Aop\x(\y_k)),\,\Aop\x(\pL(\y)) - \z(\pL(\y))\rangle \nonumber \\
\le\;& (L + L_F)\,(||\pL(\y)|| + \gamma_g)\,||\pL(\y) - \y||
        \nonumber
,\end{align}
where the first inequality uses the convexity of $g$ and $g'(\z) \in \partial g(\z)$,
the third equality uses~\eqref{eq:pd},
and the last inequality uses~\eqref{eq:pd_bound}
and Assumption~\ref{assump_g}.
\ifarxiv
        \qed
\fi
\end{proof}
\end{lemma}

Lemma~\ref{lem:pd_cost} shows
that the rate of the proximal gradient norm decrease
determines the rate of the primal-dual gap decrease of~\eqref{eq:P}
with Assumption~\ref{assump_g}.
However for problems without Assumption~\ref{assump_g},
Lemma~\ref{lem:p_d_cost} could be useful 
as an alternative measure of the convergence rate of the 
dual-based proximal gradient methods.
In addition,
the decrease of the infeasibility violation $||\Aop\x(\y) - \z(\pL(\y))||$ of~\eqref{eq:P_}
that is proportional to the proximal gradient norm decrease
based on Lemma~\ref{lem:proxcond}
could be considered for analyzing rates for such problems.

\subsection{New convergence analysis of the DPG and FDPG method}

Both DPG and FDPG have the following bound on the (dual) proximal gradient norm
\cite[Theorem 1 and Equation (5.1)]{kim:16:ala}:\footnote{
This bound is tight up to a constant for DPG~\cite{kim:16:gto}.
However, it is unknown whether or not FDPG (FISTA) 
has a bound for the proximal gradient norm decrease
that is better than the rate $O(1/k)$,
which is an interesting open question
considering that the Nesterov's fast gradient method
\cite{nesterov:83:amf}
(equivalent to FISTA for unconstrained smooth convex problems)
decreases the gradient norm with rate $O(1/k^{1.5})$ 
in~\cite{kim:16:gto}.
}
\begin{align}
||\pLkp(\y_k) - \y_k|| \le \frac{2||\y_0 - \y_*||}{k}
\label{eq:grad_fdpg}
,\end{align}
for any $L_k'$ that satisfies $\tq(\pLkp(\y_k)) \le Q_{L_k'}(\pLkp(\y_k),\y_k)$.
(Inequality~\eqref{eq:grad_fdpg}
simplifies for DPG by using $L_k' = L_{k+1}$ and $\y_{k+1} = \pLkp(\y_k)$.)
This inequality leads to new bounds on the primal-dual gap decrease
of DPG and FDPG
using Lemma~\ref{lem:pd_cost} as shown next.\footnote{
We have a primal-dual gap bound at the point $\pL(\x(\y_k))$ of FDPG
in~\eqref{eq:assumpg_cost_fdpg}
rather than that at the point $\x(\y_k)$ in~\eqref{eq:assump_cost_fdpg},
since we only know a proximal gradient norm bound 
at $\pL(\x(\y_k))$ in~\eqref{eq:grad_fdpg}.}

\begin{theorem}
Let $\{\y_k\}$ be the sequence generated 
by either DPG or FDPG.
Then with Assumption~\ref{assump_g}, 
the corresponding primal sequence defined by~\eqref{eq:primx}
satisfies
\begin{align}
&H(\x(\pLkp(\y_k))) - q(\pLkp(\y_k)) \nonumber \\
&\qquad \le (L_k' + L_F)\,(||\y_0 - \y_*|| + ||\y_*|| + \gamma_g)
        \,\frac{2||\y_0 - \y_*||}{k}
\label{eq:assumpg_cost_fdpg}
,\end{align}
for any $L_k'$ that satisfies $\tq(\pLkp(\y_k)) \le Q_{L_k'}(\pLkp(\y_k),\y_k)$.
\begin{proof}
\cite[Equation (3.6)]{beck:09:afi}
and Lemma~\ref{lem:iterbound} in Section~\ref{sec:gfdpg}
imply that the sequence $\{\y_k\}$ of both DPG and FDPG satisfy
\begin{align}
||\pL(\y_k)|| \le ||\pL(\y_k) - \y_*|| + ||\y_*||
	\le ||\y_0 - \y_*|| + ||\y_*||
\label{eq:yk_bound}
,\end{align}
where the first inequality uses the triangle inequality.
Inserting~\eqref{eq:grad_fdpg} and~\eqref{eq:yk_bound}
in Lemma~\ref{lem:pd_cost} concludes the proof.
\ifarxiv
        \qed
\fi
\end{proof}
\end{theorem}

To accelerate the rate of the primal function decrease,
the next section proposes to replace FISTA with GFPGM~\cite{kim:16:ala} 
because it decreases the proximal gradient norm with rate $O(1/k^{1.5})$.

\section{Generalized FDPG with rate $O(1/k^{1.5})$}
\label{sec:gfdpg}

The following generalized FDPG (GFDPG)
is an extension of GFPGM (with fixed $L_k$) in~\cite{kim:16:ala}
that can adopt a backtracking scheme based on~\cite{beck:09:afi}.

\begin{table}[h]
\ifarxiv
        \normalsize
\fi
\centering
\begin{tabular}{|l|}
\hline
The GFDPG method with backtracking \\
\hline
Input. Take $L_0$, $\y_0 = \w_0$, $t_0 = T_0 \in (0,\;1]$. \\
Step $k$. ($k\ge1$) \\
\quad Choose $L_k$ s.t. $L_{k-1} \le L_k$, and \\
\quad \qquad\qquad $\tq(\pLk(\w_{k-1})) \le Q_{L_k}(\pLk(\w_{k-1}),\w_{k-1})$. \\
\quad $\y_k = \pLk(\w_{k-1})$ \\
\quad Choose $t_k$ s.t. $t_k > 0$ and $t_k^2 \le T_k := \sum_{i=0}^{k}t_i$. \\
\quad $\w_k = \y_k + \frac{(T_{k-1} - t_{k-1})t_k}{t_{k-1}T_k}(\y_k - \y_{k-1})
			+ \frac{(t_{k-1}^2 - T_{k-1})t_k}{t_{k-1}T_k}(\y_k - \w_{k-1})$ \\
\hline
\end{tabular}
\end{table}

This GFDPG has the following bounds
on the dual function decrease and dual proximal gradient norm decrease
that extend
\cite[Theorems 3 and 4]{kim:16:ala}
for the GFDPG (GFPGM) with fixed $L_k$.
Note that the GFDPG and~\eqref{eq:gfdpg_cost}
reduce to FDPG and~\eqref{eq:dualcost_fdpg} respectively
when one chooses $t_k^2 = T_k$ for all $k$.

\begin{theorem}
\label{thm:bound_gfdpg}
Let $\{\y_k,\w_k\}$ be the sequence generated by GFDPG.
Then for any $k\ge1$,
\begingroup
\allowdisplaybreaks
\begin{align}
&q(\y_*) - q(\y_k) = \tq(\y_k) - \tq(\y_*) \le \frac{L_k||\y_0 - \y_*||^2}{2T_{k-1}}, 
	\label{eq:gfdpg_cost} \\
&\min\left\{\{||\y_i - \w_{i-1}||\}_{i=1}^k, ||\pLkp(\y_k) - \y_k||\right\}
\le \frac{||\y_0 - \y_*||}
        {\sqrt{\sum_{i=0}^{k-1}(T_i - t_i^2) + T_{k-1}}} 
\label{eq:gfdpg_grad}
.\end{align}
\endgroup
for any $L_k'$ that satisfies 
$\tq(\pLkp(\y_k)) \le Q_{L_k'}(\pLkp(\y_k),\y_k)$,
where $\y_i = \pLi(\w_{i-1})$.
\begin{proof}
\ifarxiv
	See Appendix~\ref{appx:a}.
        \qed
\else
	See~\ref{appx:a}.
\fi
\end{proof}
\end{theorem}

A specific version of GFDPG
requires selecting the parameters $t_k$.
We consider the choice $t_k = \frac{k+a}{a}$ for any $a>2$
that leads to the following Corollary
that provides an $O(1/k^{1.5})$ bound on the proximal gradient norm decrease
using~\cite[Corollary 2]{kim:16:ala}.

\begin{corollary}
\label{cor:gfdpga}
Let $\{\y_k,\w_k\}$ be the sequence generated by GFDPG 
with $t_k = \frac{k+a}{a}$ for any $a>2$.
Then for any $k\ge1$,
\begin{align}
&\min\left\{\{||\y_i - \w_{i-1}||\}_{i=1}^k, ||\pLkp(\y_k) - \y_k||\right\}	
	\le  \frac{a\sqrt{6}}{\sqrt{a-2}}\frac{||\y_0 - \y_*||}{k^{1.5}}
\label{eq:fdpga_grad}
\end{align}
for any $L_k'$ that satisfies 
$\tq(\pLkp(\y_k)) \le Q_{L_k'}(\pLkp(\y_k),\y_k)$,
where $\y_i = \pLi(\w_{i-1})$.
\end{corollary}

The following Lemma shows that the sequence $\{\y_k,\w_k\}$ of GFDPG
is bounded.

\begin{lemma}
\label{lem:iterbound}
Let $\{\y_k,\w_k\}$ be the sequence generated by GFDPG.
Then for any $k\ge1$,
\begin{align}
\max\{||\pLkp(\y_k)||,\,||\y_k||,\,||\w_k||\} &\le ||\y_0 - \y_*|| + ||\y_*||
\label{eq:upper}
,\end{align}
for any $L_k'$ that satisfies $\tq(\pLkp(\y_k)) \le Q_{L_k'}(\pLkp(\y_k),\y_k)$.
\begin{proof}
\ifarxiv
        See Appendix~\ref{appx:b}.
        \qed
\else
        See~\ref{appx:b}.
\fi
\end{proof}
\end{lemma}

Inserting Corollary~\ref{cor:gfdpga} and Lemma~\ref{lem:iterbound} 
to Lemmas~\ref{lem:p_d_cost} and~\ref{lem:pd_cost}
leads to the following Theorem
that bounds the primal-dual gap decrease of 
\eqref{eq:P_} and \eqref{eq:P} respectively
for GFDPG with $t_k = \frac{k+a}{a}$ for any $a>2$. 

\begin{theorem}
Let $\{\y_k,\w_k\}$ be the sequence generated by GFDPG
with $t_k = \frac{k+a}{a}$ for any $a>2$.
Then the corresponding primal sequence defined by~\eqref{eq:primx}
satisfies
\begin{align*}
&\min\left\{\begin{array}{l}
	\{\tH(\x(\y_i),\z(\y_i)) - q(\y_i)\}_{i=1}^k, \\
	H(\x(\pLkp(\y_k)),\z(\pLkp(\y_k))) - q(\pLkp(\y_k))
	\end{array}\right\} \\
&\qquad\qquad \le \frac{a\sqrt{6}}{\sqrt{a-2}}
        (L_k' + L_F)(||\y_0 - \y_*|| + ||\y_*||)
        \frac{||\y_0 - \y_*||}{k^{1.5}}
,\end{align*}
and with Assumption~\ref{assump_g} the sequence satisfies
\begin{align*}
&\min\left\{\{H(\x(\y_i)) - q(\y_i)\}_{i=1}^k,H(\x(\pLkp(\y_k))) - q(\pLkp(\y_k))\right\} \\
&\qquad\quad \le \frac{a\sqrt{6}}{\sqrt{a-2}}
        (L_k' + L_F)(||\y_0 - \y_*|| + ||\y_*|| + \gamma_g)
        \frac{||\y_0 - \y_*||}{k^{1.5}}
,\end{align*}
for any $L_k'$ that satisfies 
$\tq(\pLkp(\y_k)) \le Q_{L_k'}(\pLkp(\y_k),\y_k)$.
\end{theorem}


\begin{remark}
When one selects the total number of iterations $N$ in advance,
one can decrease the proximal gradient norm faster
than the bound~\eqref{eq:fdpga_grad}.
It is found in~\cite{kim:16:ala} that 
the following choice
\begin{align}
t_k = \begin{cases}
1, & k = 0, \\
\frac{1 + \sqrt{1 + 4t_{k-1}^2}}{2}, & k=1,\ldots,\Big\lfloor\frac{N}{2}\Big\rfloor-1, \\
\frac{N-k+1}{2}, & \text{otherwise},
\end{cases}
\end{align}
for GFPGM (and thus GFDPG)
provides the best known proximal gradient norm bound.
\end{remark}

\begin{remark}
Other accelerated proximal gradient methods
such as~\cite{ghadimi:16:agm,monteiro:13:aad}
that have $O(1/k^{1.5})$ bounds for decreasing the proximal gradient norm
could be considered instead of using GFPGM for GFDPG,
but their bounds are larger than those of GFPGM~\cite{kim:16:ala}.
\end{remark}



\section{Conclusions}
\label{sec:conc}

We provided a new analysis of the primal function decrease
of the dual-based proximal gradient methods
using the convergence analysis
of the dual proximal gradient norm.
As a consequence, we showed that using proximal gradient methods that decrease
the proximal gradient norm with rate $O(1/k^{1.5})$
leads to the same fast rate 
for the primal function (and the primal-dual gap) decrease,
improving on the previously best known rate $O(1/k)$.

\ifarxiv
	\section{Appendix}
\else
	\appendix
\fi

\ifarxiv
	\subsection{Proof of Theorem~\ref{thm:bound_gfdpg}}
	\normalsize
\else
	\section{Proof of Theorem~\ref{thm:bound_gfdpg}}
\fi
\label{appx:a}

\begin{proof}
This proof uses the fact 
that the sequence $\{\w_k\}$ of GFDPG
is equivalent to the following~\cite[Proposition 2]{kim:16:ala}: 
\begin{align}
\w_k = \frac{T_{k-1}}{T_k}\y_k
        + \frac{t_k}{T_k}\s_k
\label{eq:wk}
,\end{align}
where $\s_k := \s_{k-1} + t_{k-1}(\y_k - \w_{k-1})$.
This proof also uses the following two inequalities~\cite[Lemma 2.3]{beck:09:afi}:
\begin{align*}
&\tq(\y_{k+1}) - \tq(\y_k) 
	\le -\frac{L_{k+1}}{2}||\y_{k+1} - \w_k||^2
	- L_{k+1}\left\langle\w_k - \y_k, \y_{k+1} - \w_k\right\rangle, \\
&\tq(\y_{k+1}) - \tq(\y_*)
	\le -\frac{L_{k+1}}{2}||\y_{k+1} - \w_k||^2
	- L_{k+1}\left\langle\w_k - \y_*, \y_{k+1} - \w_k\right\rangle,
\end{align*}
and the following equality:
\begin{align*}
&||\s_{k+1} - \y_*||^2 = ||\s_k + t_k(\y_{k+1} - \w_k) - \y_*||^2 \\
=\; &||\s_k - \y_*||^2
        + 2t_k\left\langle\s_k - \y_*, \y_{k+1} - \w_k\right\rangle
        + t_k^2||\y_{k+1} - \w_k||^2
.\end{align*}

Using the above, we have
\begin{align}
&t_0(\tq(\y_1) - \tq(\y_*)) \nonumber \\
\le& -\frac{L_1t_0}{2}||\y_1 - \w_0||^2
	- L_1t_0\left\langle\w_0 - \y_*, \y_1 - \w_0\right\rangle \nonumber \\
=& -\frac{L_1}{2}\left(T_0 - t_0^2\right)||\y_1 - \w_0||^2
	+ \frac{L_1}{2}\left(||\s_0 - \y_*||^2 - ||\s_1 - \y_*||^2\right)
\label{eq:q0}
,\end{align}
and for $k\ge1$, we have
\begingroup
\allowdisplaybreaks
\begin{align}
&T_{k-1}(\tq(\y_{k+1}) - \tq(\y_k))
                + t_k(\tq(\y_{k+1}) - \tq(\y_*)) \nonumber \\
\le& -\frac{L_{k+1}T_k}{2}||\y_{k+1} - \w_k||^2 
	- L_{k+1}\left\langle T_k\w_k - T_{k-1}\y_k - t_k\y_*,
	\y_{k+1} - \w_k\right\rangle \nonumber \\
=& - \frac{L_{k+1}T_k}{2}||\y_{k+1} - \w_k||^2 
	- L_{k+1}t_k\left\langle\s_k - \y_*,\y_{k+1} - \w_k\right\rangle \nonumber \\
=& - \frac{L_{k+1}}{2}\left(T_k - t_k^2\right)||\y_{k+1} - \w_k||^2
	+ \frac{L_{k+1}}{2}\left(||\s_k-\y_*||^2 - ||\s_{k+1}-\y_*||^2\right) \nonumber 
,\end{align}
\endgroup
which becomes
\begin{align}
&\frac{1}{2}\left(T_k - t_k^2\right)||\y_{k+1} - \w_k||^2 
	+ \frac{T_k}{L_{k+1}}(\tq(\y_{k+1}) - \tq(\y_*)) \nonumber \\
\le\;& \frac{T_{k-1}}{L_{k+1}}(\tq(\y_k) - \tq(\y_*))
	+ \frac{1}{2}\left(||\s_k-\y_*||^2 - ||\s_{k+1}-\y_*||^2\right) \nonumber \\
\le\;& \frac{T_{k-1}}{L_k}(\tq(\y_k) - \tq(\y_*))
	+ \frac{1}{2}\left(||\s_k-\y_*||^2 - ||\s_{k+1}-\y_*||^2\right)
\label{eq:qk}
,\end{align}
where the last inequality uses $L_k \le L_{k+1}$.

Using a telescoping sum of~\eqref{eq:q0} and~\eqref{eq:qk}, we have 
\begin{align}
&\sum_{i=0}^{k-1}\frac{1}{2}\left(T_i - t_i^2\right)||\y_{i+1} - \w_i||^2 
	+ \frac{T_{k-1}}{L_k}(\tq(\y_k) - \tq(\y_*)) \nonumber \\
\le\;& \frac{1}{2}(||\s_0 - \y_*||^2 - ||\s_k - \y_*||^2)
\label{eq:append}
,\end{align}
which implies~\eqref{eq:gfdpg_cost}.

The condition $\tq(\pLkp(\y_k)) \le Q_{L_k'}(\pLkp(\y_k),\y_k)$
implies the following inequality~\cite[Theorem 1]{nesterov:13:gmf}:
\begin{align*}
\frac{L_k'}{2}||\pLkp(\y_k) - \y_k||^2 \le \tq(\y_k) - \tq(\pLkp(\y_k))
	\le \tq(\y_k) - \tq(\y_*)
,\end{align*}
and inserting this into~\eqref{eq:append}
leads to~\eqref{eq:gfdpg_grad}.
Note that using $0 \le \tq(\y_k) - \tq(\pLkp(\y_k))$ instead
leads to
\begin{align}
\min\{||\y_i - \w_{i-1}||\}_{i=1}^k
\le \frac{||\y_0 - \y_*||}
        {\sqrt{\sum_{i=0}^{k-1}(T_i - t_i^2)}}
\label{eq:gfdpg_grad_}
,\end{align}
which does not require computing $L_k'$ and $\pLkp(\y_k)$
unlike~\eqref{eq:gfdpg_grad},
but~\eqref{eq:gfdpg_grad_} has an upper bound
that is looser than~\eqref{eq:gfdpg_grad}.
\ifarxiv
        \qed
\fi
\end{proof}

\ifarxiv
        \subsection{Proof of Lemma~\ref{lem:iterbound}}
        \normalsize
\else
        \section{Proof of Lemma~\ref{lem:iterbound}}
\fi
\label{appx:b}

\begin{proof}
We have
\begin{align}
||\pL(\y_{k+1}) - \y_*||
&\le ||\y_{k+1} - \y_*||
\le ||\w_k - \y_*|| \nonumber \\
&\le \frac{T_{k-1}}{T_k}||\y_k - \y_*||
        + \frac{t_k}{T_k}||\s_k - \y_*|| \nonumber \\
&\le \frac{T_{k-1}}{T_k}||\y_k - \y_*||
        + \frac{t_k}{T_k}||\y_0 - \y_*|| \nonumber \\
&\le \max\{||\y_k - \y_*||,\,||\y_0 - \y_*||\}
\label{eq:yk_upper} 
,\end{align}
where the first and second inequalities use~\cite[Equation (3.6)]{beck:09:afi},
the third inequality uses~\eqref{eq:wk}, $T_k = T_{k-1} + t_k$, and a triangle inequality,
the fourth inequality uses~\eqref{eq:append},
and the last inequality uses convexity.
The inequality~\eqref{eq:yk_upper} implies
\begin{align*}
\max\{||\pL(\y_k) - \y_*||,\,||\y_k - \y_*||,\,||\w_k - \y_*||\} \le ||\y_0 - \y_*||
\end{align*}
for any $k\ge1$,
and thus inequality~\eqref{eq:upper}
follows from a triangle inequality.
\ifarxiv
        \qed
\fi
\end{proof}

\ifarxiv
\else
	\section*{Acknowledgements}
	This research was supported in part by NIH grant U01 EB018753.

	\section*{References}
\fi
\bibliographystyle{elsarticle-num} 
\bibliography{master,mastersub}

\begin{thebibliography}{10}
\expandafter\ifx\csname url\endcsname\relax
  \def\url#1{\texttt{#1}}\fi
\expandafter\ifx\csname urlprefix\endcsname\relax\def\urlprefix{URL }\fi
\expandafter\ifx\csname href\endcsname\relax
  \def\href#1#2{#2} \def\path#1{#1}\fi

\bibitem{beck:14:afd}
A.~Beck, M.~Teboulle, A fast dual proximal gradient algorithm for convex
  minimization and applications, Operations Research Letters 42~(1) (2014)
  {1--6}.

\bibitem{tseng:91:aoa}
P.~Tseng, Applications of a splitting algorithm to decomposition in convex
  programming and variational inequalities, SIAM J. Cont. Opt. 29~(1) (1991)
  {119--38}.

\bibitem{goldstein:14:fad}
T.~Goldstein, B.~O'Donoghue, S.~Setzer, R.~Baraniuk, Fast alternating direction
  optimization methods, SIAM J. Imaging Sci. 7~(3) (2014) {1588--623}.

\bibitem{beck:09:afi}
A.~Beck, M.~Teboulle, A fast iterative shrinkage-thresholding algorithm for
  linear inverse problems, SIAM J. Imaging Sci. 2~(1) (2009) {183--202}.

\bibitem{beck:09:fgb}
A.~Beck, M.~Teboulle, Fast gradient-based algorithms for constrained total
  variation image denoising and deblurring problems, IEEE Trans. Im. Proc.
  18~(11) (2009) {2419--34}.

\bibitem{pu:14:fam}
Y.~Pu, M.~N. Zeilinger, C.~N. Jones, Fast alternating minimization algorithm
  for model predictive control, in: Proc. 19th World Congress of the
  International Federation of Automatic Control, 2014, pp. {11980--6}.

\bibitem{nesterov:12:htm}
Y.~Nesterov, How to make the gradients small, Optima 88 (2012) {10--11}.

\bibitem{necoara:16:ica}
I.~Necoara, A.~Patrascu, Iteration complexity analysis of dual first order
  methods for conic convex programming, Optimization Methods and Software
  31~(3) (2016) {645--78}.

\bibitem{kim:16:ala}
D.~Kim, J.~A. Fessler, Another look at the ``{Fast Iterative
  Shrinkage/Thresholding Algorithm (FISTA)}'', arxiv 1608.03861 (2016).

\bibitem{moreau:65:ped}
J.~J. Moreau, Proximit\'e et dualit\'e dans un espace hilbertien, Bulletin de
  la Soci\'et\'e Math\'ematique de France 93 (1965) {273--99}.

\bibitem{glowinski:89}
R.~Glowinski, P.~L. Tallec, Augmented {Lagrangian} and operator-splitting
  methods in nonlinear mechanics, Soc. Indust. Appl. Math., 1989.

\bibitem{kim:16:gto}
D.~Kim, J.~A. Fessler, Generalizing the optimized gradient method for smooth
  convex minimization, arxiv 1607.06764 (2016).

\bibitem{nesterov:83:amf}
Y.~Nesterov, A method for unconstrained convex minimization problem with the
  rate of convergence {$O(1/k^2)$}, Dokl. Akad. Nauk. USSR 269~(3) (1983)
  {543--7}.

\bibitem{ghadimi:16:agm}
S.~Ghadimi, G.~Lan, Accelerated gradient methods for nonconvex nonlinear and
  stochastic programming, Mathematical Programming 156~(1) (2016) {59--99}.

\bibitem{monteiro:13:aad}
R.~D.~C. Monteiro, B.~F. Svaiter, An accelerated hybrid proximal extragradient
  method for convex optimization and its implications to second-order methods,
  SIAM J. Optim. 23~(2) (2013) {1092--1125}.

\bibitem{nesterov:13:gmf}
Y.~Nesterov, Gradient methods for minimizing composite functions, Mathematical
  Programming 140~(1) (2013) {125--61}.

\end{thebibliography}

\end{document}